\documentclass[11pt, a4paper, reqno]{amsart}

\usepackage{amssymb}
\usepackage{mathrsfs}
\usepackage{graphicx}

\usepackage[width=1.0\textwidth]{caption}

\usepackage[breaklinks]{hyperref}

\hypersetup{
unicode=true,
colorlinks=true,
linkcolor=blue,
citecolor=blue,
urlcolor=blue,
filecolor=blue,
bookmarksnumbered=true,
pdfstartview=FitH,
pdfhighlight=/N
}

\theoremstyle{plain}
\newtheorem{theorem}{Theorem}

\theoremstyle{definition}
\newtheorem{remark}{Remark}

\def\({\left(}
\def\){\right)}
\def\[{\left[}
\def\]{\right]}

\def\card{\operatorname{card}}

\def\nik{\operatornamewithlimits{Nik}}
\def\mdeg{\operatorname{deg}}
\def\mcap{\operatorname{cap}}
\def\mdiv{\operatorname{div}}
\def\supp{\operatorname{supp}}

\def\GRS{\operatorname{GRS}}

\def\NN{\mathbb N}
\def\RR{\mathbb R}
\def\CC{\mathbb C}
\def\PP{\mathbb P}

\def\HH{\mathscr H}
\def\sA{\mathscr A}
\def\RS{\mathfrak R}
\def\zz{\mathbf z}
\def\maa{\mathbf a}

\let\eps\varepsilon
\let\pfi\varphi

\let\leq\leqslant\let\geq\geqslant

\let\myh\widehat\let\myt\widetilde\let\myo\overline

\def\GRS{\operatorname{GRS}}

\begin{document}

\title{On a new approach to the problem of the zero distribution of Hermite--Pad\'e polynomials for a~Nikishin system}

\author{Sergey P. Suetin}
\address{Steklov Mathematical Institute of Russian Academy of Sciences}
\email{suetin@mi.ras.ru}
\thanks{This research was carried out with the partial support of the Russian Foundation for Basic
Research (grant no.~15-01-07531).}

\date{24.11.2017}

\begin{abstract}
A new approach to the problem of the zero distribution of Hermite--Pad\'e polynomials
of type~I for a~pair of functions $f_1,f_2$ forming a~Nikishin system is discussed.
Unlike the traditional vector approach, we give an answer in terms of a~scalar equilibrium problem
with harmonic external field, which is posed on a~two-sheeted Riemann surface.

Bibliography: \cite{Sue18} titles.

Paper: \href{http://mi.mathnet.ru/eng/tm3908}{http://mi.mathnet.ru/eng/tm3908}
\end{abstract}

\maketitle


\markright{DISTRIBUTION OF THE ZEROS OF HERMITE--PAD\'E POLYNOMIALS}

\setcounter{tocdepth}{1}
\tableofcontents


\section{Introduction and statement of the problem}\label{s1}

\subsection{}\label{s1s1}
Let
\begin{equation}
f_1(z):=\frac1{(z^2-1)^{1/2}},
\quad
f_2(z):=\int_{-1}^1\frac{h(x)}{(z-x)}\frac{dx}{\sqrt{1-x^2}},
\quad z\in D:=\myo\CC\setminus{E};
\label{1}
\end{equation}
here $E:=[-1,1]$, $h$ is a~holomorphic function on~$E$ (written $h\in\HH(E)$) of the form
$h(z)=\myh\sigma(z)$, where
\begin{equation}
\myh\sigma(z):=\int_F\frac{d\sigma(t)}{z-t},
\quad z\in\myo\CC\setminus{F},
\quad F:=\bigsqcup_{j=1}^p[c_j,d_j]\subset\RR\setminus E,
\label{2}
\end{equation}
$c_j<d_j$, $\sigma$~is a~positive Borel measure with support in~$F$
and such that $\sigma':=d\sigma/dx>0$ almost everywhere (a.e.) on~$F$.
Functions $\myh\sigma(z)$ in~\eqref{2} are called Markov functions.
Regarding the choice of branches of the function $(\,\cdot\,)^{1/2}$ and of the root $\sqrt{\,\cdot\,}$
in~\eqref{1}, see~\S\,\ref{s1s2} below.

For a tuple $[1,f_1,f_2]$ of three functions, where $f_1$ and $f_2$ are given by~\eqref{1},
and an arbitrary~$n\in\NN$, Hermite--Pad\'e polynomials of type~I $Q_{n,0},Q_{n,1},Q_{n,2}$, $\mdeg{Q_{n,j}}\leq{n}$,
$Q_{n,j}\not\equiv0$, of order~$n$ are defined (not uniquely) from the relation
\begin{equation}
R_n(z):=(Q_{n,0}\cdot 1+Q_{n,1}f_1+Q_{n,2}f_2)(z)=O\(\frac1{z^{2n+2}}\),
\quad z\to\infty.
\label{3}
\end{equation}
The purpose of the present paper is to put forward and discuss, on an example of
a~pair of functions of the form~\eqref{1}, a~new approach to
the study of the limit distribution of the zeros for Hermite--Pad\'e
polynomials of type~I as defined by~\eqref{3}.
As it is our intention to apply, in subsequent studies, this approach
to fairly general classes of analytic functions (see the result announced in~\cite{Sue18} and
Remark~\ref{rem1} below), we shall first give the notation to be used below (in this respect,
see  \cite{Sue15},~\cite{MaRaSu16},~\cite{Sue17b}).

Let $\Sigma\subset\CC$ be an arbitrary finite set, $\card\Sigma<\infty$.
We let $\sA^\circ(\Sigma)$ denote the class of all
analytic functions which are holomorphic at each point
$z_0\in\myo\CC\setminus\Sigma$, admit analytic continuation  from $z_0$
along any path~$\gamma$ in~$\myo\CC$ disjoint from~$\Sigma$,
and such that at least one point of the set~$\Sigma$ is a~branch point of this function.
For $f_1,f_2\in\sA^\circ(\Sigma)$ (under the assumption that the functions
$1,f_1,f_2$ are independent over the field~$\CC(z)$ of rational functions of~$z$ with complex coefficients),
the problem of the limit distribution of the zeros of Hermite--Pad\'e polynomials
has a~long history and in general is still unsolved
(see~\cite{Nut84}, \cite{Sta88}, \cite{Apt08},~\cite{Rak16}).
There is also no complete understanding what terms should be
employed to solve this problem.
At present, the answer to the problem of the limit distribution of the zeros of Hermite--Pad\'e polynomials
is available only for some particular classes of analytic functions
(see~\cite{GoRa81}, \cite{Nik86}, \cite{NuTr87}, \cite{GoRaSo97},
\cite{ApKuVa07}, \cite{ApLy10}, \cite{RaSu13},~\cite{MaRaSu16}).
As a~rule, the limit distribution of the zeros of Hermite--Pad\'e polynomials for a~pair
of functions $f_1,f_2$ can be described following the approach first proposed by Nuttall
(see~\cite{Nut84},~\cite{NuTr87}) in terms related to some three-sheeted
Riemann surface which in a~certain
sense\footnote{Similarly to the way
the strong asymptotics of Pad\'e polynomials
is described in terms related to the two-sheeted Riemann surface associated
(in accordance with the Stahl theory) with an arbitrary function from the class $\sA^\circ(\Sigma)$;
see \cite{Nut86}, \cite{ApYa15},~\cite{MaRaSu12}.}
is ``associated'' with the pair of functions $f_1,f_2$ (for the relation between the three-sheeted Riemann surface
with the asymptotics of Hermite--Pad\'e polynomials, see also \cite{KoKrPaSu16}, \cite{ApBoYa17},~\cite{ChKoPaSu17}.)

For a~pair of functions $f_1,f_2$ of form~\eqref{1} the above problem was solved by
Nikishin~\cite{Nik86} in 1986  (see
also \cite{Nik80}, \cite{NiSo88},~\cite{BaGeLo16}). Note that
in~\cite{Nik86} the problem was solved for an arbitrary number of functions
$f_1,f_2,\dots,f_m$ forming a~{\it Nikishin system}; a~pair of functions~\eqref{1}
is a~particular case of such a~system.
The solution of the problem of the distribution of the zeros of Hermite--Pad\'e polynomials in~\cite{Nik86}
is based on the potential theory approach developed by Gonchar and Rakhmanov~\cite{GoRa81}
in 1981 for the purposes of solving the zero distribution problem for
Hermite--Pad\'e polynomials of type~II forming an {\it Angelesco system}
(a~particular case of an arbitrary number of functions
$f_1,f_2,\dots,f_m$ was also considered in the paper~\cite{GoRa81}, in which, in particular,
the effect of {\it pushing} of the support of the equilibrium measure
inside the original orthogonality interval was discovered; see
also~\cite{Rak18}).
Within the framework of this vector approach, the answer for a~pair of functions~\eqref{1}
is given in terms of a~vector-equilibrium measure $\vec{\lambda}=(\lambda_1,\lambda_2)$
supported on the vector-compact set
$(E,F)$ (that is,
$\supp{\lambda_1}\subset{E},\supp{\lambda_2}\subset{F}$).
The equilibrium conditions are determined by the interaction matrix of measures
$M_{\nik}=
\begin{pmatrix} 2&-1\\-1&2\end{pmatrix}$, which is known as the {\it Nikishin matrix}.
The solution of the problem is a~unique
vector-measure $\vec{\lambda}=(\lambda_1,\lambda_2)$ with support on the
 vector-compact set  $(E,F)$; this measure  is {\it extremal}
for the energy functional defined by the logarithmic kernel and the interaction matrix $M_{\nik}$
(see~\cite{NiSo88},~\cite{GoRaSo97},~\cite{Apt08},~\cite{Lap15}).
The extremal vector-measure
$\vec{\lambda}=(\lambda_1,\lambda_2)$ is completely characterized by the
{equilibrium} condition for the corresponding vector potential and the vector-compact set
$(E,F)$; for more details, see \cite{Apt08},~\cite{Lap15}.

Note that, for arbitrary functions $f_1,f_2\in\sA^\circ(\Sigma)$, the
problem of the limit distribution of the zeros of the corresponding
Hermite--Pad\'e polynomials turns out to be equivalent to the problem
of the limit distribution of the zeros of polynomials satisfying some
{\it non-Her\-mi\-tian} orthogonality conditions
(see~\cite{GoRa87},~\cite{Rak12}, \cite{Rak16},~\cite{Rak18}). The
characteristic feature of non-Hermitian orthogonality conditions is
that the contour of integration is {\it not fixed a~priori}, but rather
lies in some class of ``admissible'' contours. The following heuristic
conclusion can be made based on a~series of particular cases
investigated so far: in this class, there exists a~{\it unique
``optimal''} contour\footnote{Here and below, by a~contour we shall
mean a~composite contour consisting of a~finite number of closed curves
and splitting the Riemann sphere into a~finite number of domains; see
\cite{Bus13}, \cite{Bus15}.} attracting in the limit the zeros of
Hermite--Pad\'e polynomials. The  ``optimality'' property of a~contour
is formulated in terms of the corresponding vector equilibrium problems
of potential theory. This optimal contour possesses a~certain vector
$S$-property, which completely characterizes it in the class of
admissible vector-contours. In modern terms, such a~contour is called
an $S$-curve or an  $S$-compact set (see~\cite{Rak12}).

The concept of an $S$-compact set was first introduced by  H.~Stahl in the 1985--1986s
(see~\cite{Sta97b} and~\cite{Sta12} and there references given therein)
when considering the problem of the limit distribution of the zeros and poles of Pad\'e approximants
 in the class of multivalued analytic functions
$\sA^\circ(\Sigma)$.
In~1987 Gonchar and Rakhmanov~\cite{GoRa87}, in their solution of the ``$1/9$ conjecture'',
developed a~different approach to the problem of the limit distribution of the zeros
of non-Hermitian orthogonal polynomials. This approach is based on the
scalar equilibrium problem, but with the so-called ``\textit{external field}''
defined by a~harmonic function
(more general external fields and the corresponding $S$-curves were considered in~\cite{Rak12}).
This new approach was used in 2012--2015 by Buslaev
\cite{BuMaSu12}--\cite{Bus15} to solve the problem of the limit distribution of the zeros and poles
of multivalued Pad\'e approximants.
Here, the potential of a~negative unit charge concentrated at a~finite number
of interpolation nodes appears naturally
as an external field (see also \cite{Bus15b}, \cite{BuSu16},~\cite{Bus16}).

The class of methods developed by H.~Stahl, A.~A. Gonchar, and
E.~A. Rakh\-ma\-nov in the 1980s for the purpose of studying the limit distribution of the zeros
of non-Hermitian orthogonal polynomials is called at present the Gonchar--\allowbreak Rakhmanov--\allowbreak Stahl method (or briefly
the  $\GRS$-{\it method});
see~\cite{Sue15}, \cite{MaRaSu16}, \cite{Rak16}, \cite{Rak18}.

The purpose of the present paper is, by using an example of two functions $f_1$ and $f_2$ of the form~\eqref{1},
put forward and discuss a~new approach to the
problem of the limit distribution of the zeros of Hermite--Pad\'e polynomials, which in a~certain sense further develops
the approach of A.~A.~Gonchar and E.~A.~Rakh\-ma\-nov employed in their solution
of the ``$1/9$ conjecture''. Namely, the limit
distribution of the zeros of the polynomial $Q_{n,2}$ as $n\to\infty$
will be characterized in terms related to some scalar
potential theory equilibrium problem (but with external field), which in addition
is posed not on the Riemann sphere $\myo\CC$, but rather on the two-sheeted Riemann surface of the function $w^2=z^2-1$.
This is the principal distinguishing feature of the approach of the present paper
from the standard method based on the vector equilibrium problem posed on the Riemann sphere.


Let us clarify the choice of the pair of functions \eqref{1}  to illustrate the new approach and
the fact that here we speak only about the distribution of the zeros of the polynomial
$Q_{n,2}$.

The thing is, on the one hand, as we have already mentioned, in the class $\sA^\circ(\Sigma)$ the
problem of the distribution of the zeros of Hermite--Pad\'e polynomials for an
arbitrary pair of independent functions $f_1,f_2\in\sA^\circ(\Sigma)$ is not yet solved
and it is even unclear what terms should be employed to find its solution
(for conjectures in this direction, see \cite{Nut84},~\cite{Apt08},
\cite{Sta88},~\cite{Rak16}). In particular, there is no solution in this problem
even for a~pair of functions with two branch points, of which each is in
``the general position''. On the other hand, for the Pad\'e polynomials
$P_{n,0},P_{n,1}$, $\mdeg{P_{n,j}}\leq{n}$, $P_{n,j}\not\equiv0$,
as defined from the relations
\begin{equation}
(P_{n,0}+P_{n,1}f)(z)=O\(\frac1{z^{n+1}}\),
\quad z\to\infty,
\label{p1}
\end{equation}
where $f\in\HH(\infty)$,
Stahl's theory is valid \footnote{Stahl's theory is much more general and
can be applied to any multivalued analytic function, whose singular set is
of zero logarithmic capacity.} for an arbitrary function~$f$ from the class
$\sA^\circ(\Sigma)$. This leads to the following fairly natural argument:
one of the functions, $f_1$ say, in the relation~\eqref{3} defining the Hermite--Pad\'e polynomials,
should be taken as simple as possible (retaining the independence of two functions
$f_1,f_2$) with the aim at maximally extending the calculus to a~larger class of
functions that contains the second function $f_2$.
In this way, using the approach proposed here, we managed to substantially enlarge
the class of functions containing the function~$h$ in representation~\eqref{1}.
Namely, for an {\it arbitrary} function
$h\in\sA^\circ(\Sigma)$, where $\Sigma\subset\myo\CC\setminus{E}$, it is
possible
to characterize completely the problem of the limit distribution of the zeros of the polynomials $Q_{n,2}$
in terms of the same scalar potential theory equilibrium problem  with an
external field.
This problem is posed on the same two-sheeted Riemann surface of the function $w^2=z^2-1$
in a~similar way as in the present paper. The principal difference is that in the general
case it is first required to establish the existence of an appropriate
$S$-compact set~$F$ corresponding to the problem under consideration and which replaces
the union of a~finite number of closed intervals
(see~\eqref{2}). The corresponding result was announced in~\cite{Sue18};
the author intends to give the proof of
this  result in a~separate paper.

We note the papers \cite{RaSu13},~\cite{Sue15} and~\cite{KoSu15},
in which
the equilibrium problem for a~mixed Green-logarithmic potential
was employed
for the study of the limit distribution of the zeros
of Hermite--Pad\'e polynomials for a~tuple
$[1,f_1,f_2]$,
where a~pair of functions $f_1,f_2$ forms a~generalized (complex) Nikishin system
(see also~\cite{BuSu15},
\cite{Sue16},~\cite{MaRaSu16},~\cite{Rak16}).
The method of investigation proposed in the present paper is different
from that of \cite{RaSu13},~\cite{Sue15} and~\cite{KoSu15}.
Some precursor considerations and results that eventually culminated in the statement of
the potential theory equilibrium problem
on the Riemann surface $w^2=z^2-1$ were obtained by the author in~\cite{Sue17}.

It also should be mentioned about the papers by H.~Stahl with coauthors \cite{BaStYa12} and~\cite{LuSiSt15},
in which potentials pretty close to those used in the present paper were used.
However, as far as the author is aware, such potentials have not been applied before
in the study of the distribution of the zeros of Hermite--Pad\'e polynomials.

It is worth pointing out that in the present paper we
discuss and examine only the case of
``diagonal'' (that is, of the same degree) Hermite--Pad\'e polynomials
of type~I. The nondiagonal case, as well as in the case of Hermite--Pad\'e polynomials of type~II
merits special consideration within the framework of the new approach proposed here
(of course, if such a~research will prove feasible).

\begin{center}
\begin{figure}[!ht]
\includegraphics[scale=1.6]{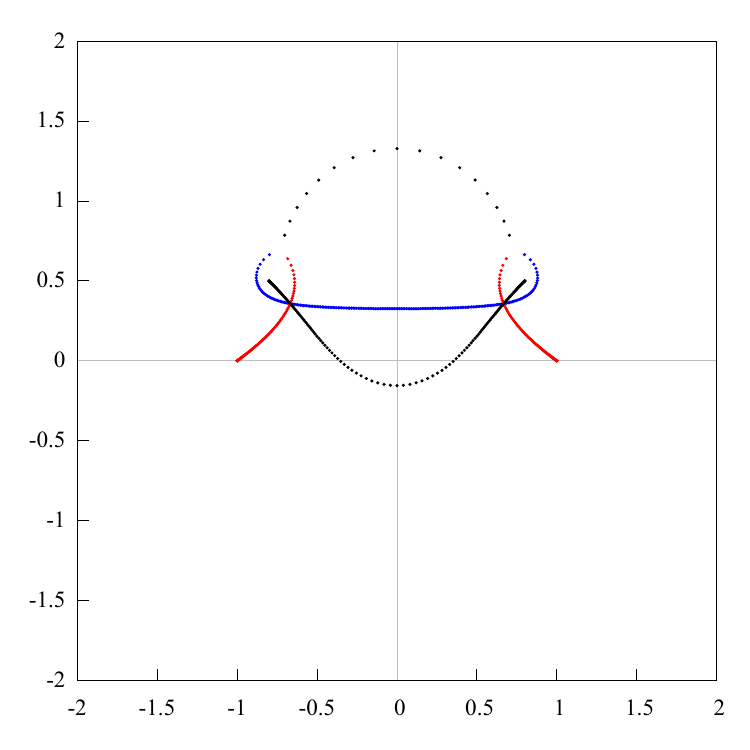}
\caption{Zeros of diagonal Hermite--Pad\'e polynomials of type~I
$Q_{200,0}$ (blue points),
$Q_{200,1}$ (red points),
$Q_{200,2}$ (black points)
for the tuple of functions $[1,f_1,f_2]$, where
$f_1(z):=(z^2-1)^{-1/2}$,
$f_2(z):=\bigl((z-.8-.5i)(z+.8-.5i)\bigr)^{-1/2}$,
forming an Angelesco system.
No theoretical justification of such behavior of the zeros of
Hermite--Pad\'e polynomials of type~I has not yet been found to date.
}
\label{Fig_hp1}
\end{figure}
\end{center}

\begin{center}
\begin{figure}[!ht]
\includegraphics[scale=1.6]{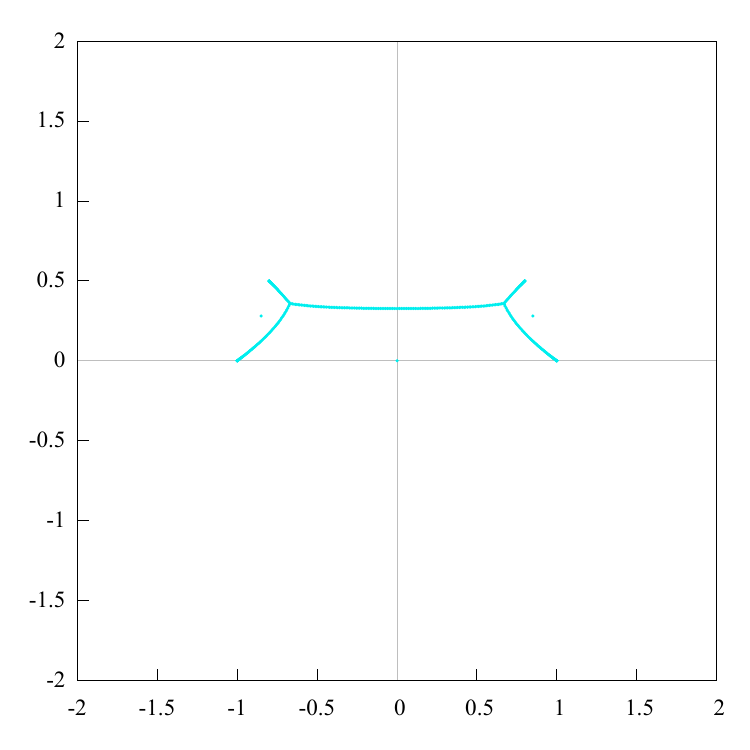}
\caption{Zeros of the denominator of diagonal Hermite--Pad\'e approximants of type~II
$P_{400}$ (light blue points)
for the tuple of functions $[1,f_1,f_2]$, where
$f_1(z):=(z^2-1)^{-1/2}$,
$f_2(z):=\bigl((z-.8-.5i)(z+.8-.5i)\bigr)^{-1/2}$,
forming an Angelesco system.
Theoretical justification of such behavior of zeros of
Hermite--Pad\'e polynomials of type~II was obtained in~\cite{ApKuVa07}.
}
\label{Fig_hp2}
\end{figure}
\end{center}

The fact that the problems on the
distribution of the zeros                         of Hermite--Pad\'e polynomials of type~I and type~II
are substantially different and in general call for different approaches and methods of
investigation is well illustrated in Figs.~\ref{Fig_hp1}--\ref{Fig_hp2}, which were derived for the
pair of functions
$$
f_1(z):=\frac1{(z^2-1)^{1/2}},\quad
f_2(z):=\frac1{\bigl((z-.8-.5i)(z+.8-.5i)\bigr)^{1/2}},
$$
forming an Angelesco  system.

\smallskip

The author is grateful to the referee for the many helpful comments and suggestions
which led to a great improvement in the presentation of the paper and for calling
his attention to the papers \cite{Apt99} and~\cite{LoVa16}.

\subsection{}\label{s1s2}
We shall require the following notation and definitions.
We set  $D:=\myo\CC\setminus{E}$,
\begin{equation}
\pfi(z):=z+(z^2-1)^{1/2},\quad z\in D,
\label{5}
\end{equation}
where we choose the branch of the root function such that $(z^2-1)^{1/2}/z\to1$ as $z\to\infty$.
For $x\in (-1,1)$, by $\sqrt{1-x^2}$ we shall understand the
positive square root: $\sqrt{b^2}=b$ for $b\geq0$.

Given an arbitrary polynomial $Q\in\CC[z]$, $Q\not\equiv0$, by
\begin{equation*}
\chi(Q):=\sum_{\zeta:Q(\zeta)=0}\delta_\zeta
\end{equation*}
we shall mean the counting measure of the zeros of the polynomial~$Q$ (counting multiplicities).
In what follows, given an arbitrary $n\in\NN$, we denote by $\PP_n:=\CC_n[z]$
the class of all algebraic polynomials of degree $\leq{n}$ with complex coefficients.

We let $\RS_2$ denote the two-sheeted Riemann surface of the function $w^2=z^2-1$
regarded as a two-sheeted covering of the extended complex plane~$\myo\CC$
with branch points at $z=\pm1$. Each (open) sheet of the Riemann surface $\RS_2$ is the
Riemann sphere cut along the interval~$E$, the opposite sides of cuts from different sheets
being identified. The first
(open) sheet $\RS^{(1)}$ of the Riemann surface $\RS_2$ is that on which
$w=(z^2-1)^{1/2}\sim z$ as $z\to\infty$; on the second sheet $\RS^{(2)}$
$w=-(z^2-1)^{1/2}\sim -z$ as $z\to\infty$. A~point~$\zz$ on the Riemann surface~$\RS_2$
is the pair $(z,w)=\zz\in\RS_2$. The canonical projection
$\pi$, $\pi\colon\RS_2\to\myo\CC$, is defined in the standard way:
$\pi(\zz)=z$. Note that  $\RS_2$ is a~Riemann surface of zero genus, and so
any divisor~$d$ of degree~$0$ on~$\RS_2$ is a~principal one; that is,
there exists a~meromorphic function on~$\RS_2$ whose divisor of the zeros and poles
coincides with~$d$. From a~given divisor of degree~$0$ such a~meromorphic
function is defined uniquely up to a~nontrivial multiplicative
constant (for a~more detailed account of these and other aspects of Riemann surfaces, see~\cite{Chi06}).

Thus, the function $z+w$, which is meromorphic on the Riemann surface $\RS_2$, will be denoted by
$\Phi(\zz):=z+w$. Points of the Riemann surface lying on the first (open) sheet $\RS^{(1)}$
will be denoted by~$z^{(1)}$; by $z^{(2)}$ we denote points from the second sheet $\RS^{(2)}$.
So, $z^{(1)}=(z,(z^2-1)^{1/2})$,
$z^{(2)}=(z,-(z^2-1)^{1/2})$, $\pi(\RS^{(1)})=\pi(\RS^{(2)})=D$.

The following identity\footnote{This identity, for
all its undoubted simplicity, was first used very effectively
in~\cite{GoSu04}, formulas (15),~(67).}
is easily verified for $\zz,\maa\in\RS_2\setminus\Gamma$
\begin{equation}
z-a\equiv -\frac{[\Phi(\zz)-\Phi(\maa)][1-\Phi(\zz)\Phi(\maa)]}
{2\Phi(\zz)\Phi(\maa)}.
\label{4}
\end{equation}
Indeed, each of the functions on the right and left of~\eqref{4}
is meromorphic on the Riemann surface $\RS_2$. The divisor $z-a$ of the left-hand side
can be easily evaluated to be equal to $d_1=-\infty^{(1)}-\infty^{(2)}+a^{(1)}+a^{(2)}$.
For the divisor of the right-hand side, we also have
$d_2=-\infty^{(1)}-\infty^{(2)}+a^{(1)}+a^{(2)}$. Hence, these two
functions are identically equal except for a~multiplicative constant, which can be
easily calculated.

From~\eqref{4} we have, in particular, the identity
\begin{equation}
z-a\equiv -\frac{[\pfi(z)-\pfi(a)][1-\pfi(z)\pfi(a)]}
{2\pfi(z)\pfi(a)},
\label{4.2}
\end{equation}
which holds for $z,a\in D$. The following identity
\begin{equation}
\Phi(z^{(1)})\Phi(z^{(2)})\equiv1,\quad z\in D,
\label{4.3}
\end{equation}
can also be easily verified.

Let $M_1(F)$ be the space of all unit positive Borel measures supported
on a~compact set~$F$. Given an arbitrary measure $\mu\in M_1(F)$, we define by
\begin{equation}
V^\mu(z):=\int_F\log\frac1{|z-t|}d\mu(t)
\label{7}
\end{equation}
the logarithmic potential of~$\mu$,
\begin{equation}
I(\mu):=\iint_{F\times F}\log\frac1{|z-t|}d\mu(z)\,d\mu(t)
=\int_F V^\mu(z)\,d\mu(z)
\label{8}
\end{equation}
is the corresponding energy functional. By $M_1^\circ(F)\subset M_1(F)$ we shall
denote the space of measures with finite energy, $I(\mu)<\infty$.
We recall the positivity property of logarithmic energy\footnote{More precisely, the
positivity of the logarithmic kernel.}
with respect to neutral charges:
\begin{equation}
I(\mu-\nu)\geq0\quad\forall\mu,\nu\in M_1^\circ(F)
\quad\text{and}\quad I(\mu-\nu)=0 \Leftrightarrow\mu=\nu.
\label{9}
\end{equation}
For an account of these and other properties of logarithmic potentials employed in the present paper, see~\cite{Lan66}.

For a~measure $\mu\in M_1(F)$, we set\footnote{It is clear that $\psi(z)=\log{\pfi(z)}$
for $z\in\RR\setminus{E}$.}
\begin{equation}
P^\mu(z):=
\int_F\log\frac{|1-\pfi(z)\pfi(t)|}{|z-t|^2}\,d\mu(t),
\quad \psi(z):=\log|\pfi(z)|,
\label{10}
\end{equation}
and define
\begin{equation}
\begin{aligned}
J(\mu):&=
\iint_{F\times F}\log\frac{|1-\pfi(z)\pfi(t)|}{|z-t|^2}\,d\mu(z)\,d\mu(t)\\
&=\int_F P^\mu(z)\,d\mu(z),\\
J_\psi(\mu):&=
\iint_{F\times F}\biggl\{
\log\frac{|1-\pfi(z)\pfi(t)|}{|z-t|^2}
+\psi(z)+\psi(t)\biggr\}\,d\mu(z)\,d\mu(t)\\
&=\int_F P^\mu(z)\,d\mu(z)+2\int_F\psi(z)\,d\mu(z).
\end{aligned}
\label{11}
\end{equation}
From identity \eqref{4.2} we have the following equality, which holds for $z,\zeta\in D$,
\begin{equation}
\log\frac{|1-\pfi(z)\pfi(\zeta)|}{|z-\zeta|^2}=
\log\frac1{|z-\zeta|}+\log\frac1{|\pfi(z)-\pfi(\zeta)|}
+\log2+\psi(z)+\psi(\zeta).
\label{12}
\end{equation}
Potentials with kernels of the form
$$
\log\frac1{|z-\zeta|}+\log\frac1{|v(z)-v(\zeta)|},
$$
where $z,\zeta\in[A,B]\subset\RR$, $v(z)$ is an arbitrary nondecreasing function
on $[A,B]$, were considered in the paper~\cite{LuSiSt15}, however, the author of
the present paper is unaware of any applications of such potentials in the theory of Hermite--Pad\'e polynomials.

\subsection{}\label{s1s3}
The main results of the present paper are Theorems~\ref{th1} and~\ref{th2}.

\begin{theorem}\label{th1}
In the class $M_1^\circ(F)$, there exists a~unique measure
$\lambda=\lambda_{F}\in M_1^\circ(F)$ such that
\begin{equation}
J_\psi(\lambda)=\min_{\mu\in M_1(F)}J_\psi(\mu).
\label{13}
\end{equation}
The measure $\lambda$ is completely characterized by the following equilibrium condition:
\begin{equation}
P^\lambda(z)+\psi(z)
\begin{matrix}
\,\equiv w_F, & z\in S(\lambda),\\
\,\geq w_F, & z\in F\setminus S(\lambda).
\end{matrix}
\label{14}
\end{equation}
\end{theorem}

\begin{theorem}\label{th2}
Let $f_1$ and  $f_2$ be functions given by the representations~\eqref{1} and let $Q_{n,2}$
be the Hermite--Pad\'e polynomial defined by~\eqref{3}. Then
\begin{equation}
\frac1n\chi(Q_{n,2})\to\lambda,\quad n\to\infty.
\label{15}
\end{equation}
\end{theorem}

The convergence in \eqref{15} shall be understood in the sense of weak convergence in the space of
measures.
It may be pointed out once more that the assertion of Theorem~\ref{th2}
on the existence of the limit distribution of the zeros of the polynomials $Q_{n,2}$
is not new (see, first of all,~\cite{Nik86}, and also~\cite{GoRaSo97},~\cite{ApLy10}).
The new point here is the characterization of this limit distribution in terms
the scalar equilibrium problem \eqref{13}--\eqref{14}. This was achieved by posing the corresponding
potential theory problem not on the Riemann sphere, but on the two-sheeted
Riemann surface of the function $w^2=z^2-1$.

\section{Proof of Theorem~\ref{th1}}\label{s2}

\subsection{}\label{s2s1}
Let $U\supset F$ be some neighborhood of the compact set~$F$ such that
$U\cap{E}=\varnothing$. For all $\mu\in M_1(F)$, the function
$$
\int_F\log|1-\pfi(z)\pfi(t)|\,d\mu(t)
$$
is harmonic in $U$ and  the potential  $V^\mu(z)$ is a~superharmonic function
in~$U$, and hence, since $\mcap{F}>0$, $M_1(F)$
is compact in the weak topology, and using the principle of descent for logarithmic potentials
(see~\cite{Lan66}, Ch.~I, \S\,3, Theorem~1.3), we see that there exists a~measure $\lambda\in M_1^\circ(F)$
satisfying equality~\eqref{13}. Using identity~\eqref{12}, one can easily prove the convexity
of the energy functional $J_\psi(\cdot)$,
\begin{equation}
J_\psi\(\frac{\mu+\nu}2\)\leq\frac12\bigl[J_\psi(\mu)+J_\psi(\nu)\bigr]
\quad \forall \mu,\nu\in M_1(F);
\label{18}
\end{equation}
moreover,
\begin{gather}
J(\mu-\nu)=2J_\psi(\mu)+2J_\psi(\nu)-4J_\psi\(\frac{\mu+\nu}2\),
\label{19}\\
J(\mu-\nu)=\iint\limits_{F\times F}
\biggl\{\log\frac1{|z-\zeta|}+\log\frac1{|\pfi(z)-\pfi(\zeta)|}
\biggr\}\,d(\mu-\nu)(z)\,d(\mu-\nu)(\zeta),
\label{20}
\end{gather}
for all  $\mu,\nu\in M_1^\circ(F)$. As a~direct corollary of~\eqref{18}--\eqref{20}
we see that the functional
 $J(\cdot)$ is positive on neutral charges (cf.~\eqref{9}),
$$
J(\mu-\nu)\geq0\quad \forall\mu,\nu\in M_1^\circ(F)
\quad\text{and}\quad J(\mu-\nu)=0\Leftrightarrow\mu=\nu.
$$
Furthermore, the following equalities are easily verified:
\begin{gather}
J(\mu)=\iint_{F\times F}\biggl\{\log\frac1{|z-\zeta|}+
\log\frac1{|\pfi(z)-\pfi(\zeta)|}\biggr\}+\log{2}+2\int\pfi(z)\,d\mu(z),
\notag\\
J_\psi(\mu)=\iint_{F\times F}\biggl\{\log\frac1{|z-\zeta|}+
\log\frac1{|\pfi(z)-\pfi(\zeta)|}\biggr\}+\log{2}+4\int\pfi(z)\,d\mu(z).
\end{gather}

\subsection{}\label{s2s2}
Arguing as in Lemma~6 of~\cite{GoRa87}, we can now prove the equilibrium property~\eqref{14} of the extremal measure~$\lambda$
by using the above equalities and the positivity property of the functional $J(\cdot)$

Indeed\footnote{For completeness of presentation, we give
the proof of~\eqref{14}, cf.~Lemma~6 of~\cite{GoRa87}.}
one verifies directly that
\begin{equation}
J_\psi(\eps\nu+(1-\eps)\lambda)-J_\psi(\lambda)
=2\eps\int_F(P^\lambda+\psi)(z)\,d(\nu-\lambda)+\eps^2J(\nu-\lambda)
\label{7.1}
\end{equation}
for any $\eps>0$ and measure  $\nu\in M_1^\circ(F)$. It follows that the
minimizing measure~$\lambda$ is the only measure from $M_1^\circ(F)$ satisfying the condition
\begin{equation}
\int_F(P^\lambda+\psi)\,d(\nu-\lambda)\geq0\quad
\forall \nu\in M_1^\circ(F).
\label{7.2}
\end{equation}
In the actual fact, \eqref{7.2} is an immediate consequence of \eqref{18},~\eqref{9} and~\eqref{7.1} as $\eps\to0$.
On the other hand, since the energy functional $J(\cdot)$ is positive on neutral charges,
we have $J(\nu-\lambda)\geq0$ for any measure $\nu\in
M_1^\circ(F)$. An appeal to~\eqref{7.1} with $\eps=1$ shows that any measure $\nu\in
M^\circ_1(F)$ satisfying~\eqref{7.2} minimizes the energy integral $J_\psi(\cdot)$.
If a~measure $\lambda$ satisfies condition~\eqref{7.2}, then it obeys the
equilibrium relations~\eqref{14} with
$$
w_F:=\int_F(P^\lambda+\psi)\,d\lambda.
$$
Indeed, if $P^\lambda(x)+\psi(x)<w_F$ on a~closed set $e\subset
F$, $\mcap(e)>0$, then there exists $\nu\in M_1^\circ(e)$, for which
$\displaystyle\int_F(P^\lambda+\psi)(x)\,d\nu(x)<w_F$, which
shows that \eqref{7.2}~is violated. Hence $(P^\lambda+\psi)(x)\geq w_F$
everywhere on the (regular) compact set~$F$.
If $(P^\lambda+\psi)(x)>w_F$ on a~nonempty set $e\subset S(\mu)$,
then the inequality
$\displaystyle\int(P^\lambda+\psi)(x)\,d\lambda(x)>w_F$ is secured by the
lower semi-continuity of the function $(P^\lambda+\psi)(z)$, contradicting the definition
of~$w_F$.

If $\lambda$ is an equilibrium measure, then $P^\lambda+\psi\leq w_F$ everywhere on
$S(\lambda)$, which shows that $\lambda\in M_1^\circ(F)$.
Since the sets of zero inner capacity play no role in integration with
respect to measures in~$M_1^\circ(F)$, we obtain~\eqref{7.2}. Finally,
$(P^\lambda+\psi)(z)\equiv w_F$ on~$S(\lambda)$, because
$F$~is a~regular compact set.

Thus, the extremal measure $\lambda$, and only this measure, satisfies
the equilibrium conditions~\eqref{14}. This proves Theorem~1.

\smallskip

Note that $F$ is a~regular compact set, and hence the equilibrium measure
is characterized by the equality
$$
\min_{z\in F}(P^\lambda+\psi)(z)=\max_{\mu\in M_1(F)}
\min_{z\in F}(P^\mu+\psi)(z).
$$

\section{Proof of Theorem \ref{th2}}\label{s3}

\subsection{}\label{s3s1}
From \eqref{3} we have the relation
\begin{equation}
0=\int_\gamma(Q_{n,0}+Q_{n,1}f_1+Q_{n,2}f_2)(z)q(z)\,dz=
\int_\gamma(Q_{n,1}f_1+Q_{n,2}f_2)(z)q(z)\,dz,
\label{22}
\end{equation}
which holds for any polynomial $q\in\PP_{2n}$; in~\eqref{22}\enskip  $\gamma$~is an
arbitrary contour separating the interval~$E$ from the infinity point
$z=\infty$.

Let $P_{n,}$ and $P_{n,1}$ be the Pad\'e polynomials for the function $f_1$; that is,
$\mdeg{P_{n,j}}\leq \nobreak n$, $P_{n,j}\not\equiv0$, and
\begin{equation}
H_n(z):=(P_{n,0}+P_{n,1}f_1)(z)=O\(\frac1{z^{n+1}}\),\quad z\to\infty.
\label{23}
\end{equation}
It is known that $P_{n,1}=T_n$ are Chebyshev polynomials of the first kind
that are orthogonal on the interval~$E$ with the weight $1/\sqrt{1-x^2}$, $H_n$
is the corresponding function of the second kind. We shall assume that the Chebyshev polynomials are
normalized as follows: $T_n(z)=2^nz^n+\dotsb$. Hence, for the functions of the
second kind $H_n$, we have
\begin{equation}
H_n(z)=\frac{\varkappa_n\pfi'(z)}{\pfi^{n+1}(z)},\quad\varkappa_n\neq0,
\quad
H_n(z)=\frac1{2\pi i}\int_E\frac{T_n(x)\Delta f_1(x)}{x-z}\,dx,\quad
z\in D,\label{24}
\end{equation}
\begin{align}
\Delta H_n(x):&=H_n(x+i0)-H_n(x-i0)\notag\\
&=T_n(x)\Delta f_1(x)=T_n(x)\frac2{i\sqrt{1-x^2}},\quad x\in(-1,1).
\label{25}
\end{align}
Besides, the polynomials $T_n$ and the functions of the second kind $H_n$ satisfy the
same second-order recurrence relation, but with different initial data
\begin{equation}
y_k=2zy_{k-1}-y_{k-2}, \quad k=1,2,\dots,
\label{26}
\end{equation}
where one should put $y_{-1}\equiv0$, $y_0\equiv1$ for the polynomials $T_k$ and
$y_{-1}\equiv1$, $y_0=f_1(z)=1/(z^2-1)^{1/2}$ for the functions of the second
kind $H_k$.
We have
$$
\int_\gamma p(z)f_1(z)T_{n+j}(z)\,dz=0,\quad j=1,2,\dots,n
$$
for any polynomial $p\in\PP_n$, and so from \eqref{22} with
$q=T_{n+1},\dots,T_{2n}$ it follows that
\begin{equation}
\int_\gamma Q_{n,2}(z)f_2(z)T_{n+j}(z)\,dz=0,\quad
j=1,2,\dots,n.
\label{27}
\end{equation}
Next, using \eqref{27} and the definition~\eqref{1} of the function $f_2$, we have
\begin{equation}
\int_E Q_{n,2}(x)T_{n+j}(x)\frac1{\sqrt{1-x^2}}h(x)\,dx=0,
\quad j=1,\dots,n.
\label{28}
\end{equation}
In view of~\eqref{25}, the above relation is equivalent to the relation
\begin{equation}\label{eq31x}
\int_\gamma Q_{n,2}(z)H_{n+j}(z)h(z)\,dz=0,\quad
j=1,\dots,n,
\end{equation}
where $\gamma$ is an arbitrary contour separating the interval~$E$ from the compact set~$F$.
Since $h(z)=\myh\sigma(z)$, relation \eqref{eq31x} can be easily written in the form
\begin{equation}
\int_F Q_{n,2}(x)H_{n+j}(x)\,d\sigma(x)=0,
\quad j=1,\dots,n.
\label{29}
\end{equation}
These orthogonality relations\footnote{In view of the representation
$H_n(z)=\pfi'(z)/\pfi^{n+1}(z)$, the orthogonality relations~\eqref{29}
are similar to those considered in~\cite{LuSiSt15}.}
will play a key role in the subsequent analysis of
the limit distribution of the zeros of the polynomials $Q_{n,2}$.

Let $N$, $0\leq N\leq n$, be an arbitrary natural number. We shall assume
without loss of generality that $N=2m$ is an even number (the case of an odd~$N$
is treated similarly). Given arbitrary complex numbers
$c_1,\dots,c_N\in\CC$, consider the sum
$$
\sum_{j=1}^N c_jH_{n+j}(z).
$$
By using the recurrence relations \eqref{26}, this sum can be easily written as
\begin{equation}
\sum_{j=1}^Nc_jH_{n+j}(z)=q_{m,1}(z)H_{n+m+1}(z)+q_{m,2}(z)H_{n+m}(z),
\label{30}
\end{equation}
where $q_{m,1},q_{m,2}\in\PP_{m-1}$ are polynomials of degree $\leq{m-1}$. Since
the constants $c_1,\dots,c_N$ in~\eqref{30} are arbitrary,
it is easily verified that the polynomials $q_{m,1}$ and $q_{m,2}$
can also be chosen arbitrarily. So, using~\eqref{30}, relations~\eqref{29} can be
written in the following equivalent form
\begin{equation}
\int_F Q_{n,2}(x)\bigl\{q_{m,1}(x)H_{n+m+1}(x)+q_{m,2}(x)H_{n+m}(x)\bigr\}
\,d\sigma(x)=0
\label{31}
\end{equation}
with arbitrary polynomials $q_{m,1}\in\PP_{m-1}$ and $q_{m,2}\in\PP_{m-1}$.
Now,
from~\eqref{31} and the available properties of the functions of the second
kind (see~\eqref{24}), we have
\begin{align}
0&=\int_F Q_{n,2}(x)\biggl\{q_{m,1}(x)\frac{H_{n+m+1}}{H_{n+m}}(x)
+q_{m,2}(x)\biggr\}H_{n+m}(x)\,d\sigma(x)\notag\\
&=\int_F Q_{n,2}(x)\biggl\{q_{m,1}(x)
\frac{\varkappa_{n+m+1}}{\varkappa_{n+m}\pfi(x)}+q_{m,2}(x)\biggr\}
\frac{\varkappa_{n+m}\pfi'(x)}{\pfi^{n+m+1}(x)}\,d\sigma(x).
\label{32}
\end{align}



Now, using the definition of the function $\Phi(\zz)$ (see sec.~\ref{s1s2}),
which is meromorphic on the Riemann surface $\RS_2$, we get the following orthogonality relation
\begin{equation}
\int_F Q_{n,2}(x)\Bigl\{q_{m,1}(x)\Phi(x^{(2)})+q_{m,2}(x)\Bigr\}
\pfi'(x)\Phi(x^{(2)})^{n+m+1}\,d\sigma(x)=0,
\label{33}
\end{equation}
which holds for any polynomials $q_{m,1},q_{m,2}\in\PP_{m-1}$.

\subsection{}\label{s3s2}
We now set
\begin{equation}
g_n(\zz):=q_{m,1}(z)\Phi(\zz)+q_{m,2}(z),
\label{34}
\end{equation}
where it is assumed that $\mdeg{q_{m,1}}=\mdeg{q_{m,2}}=m-1$.
Then, for the divisor of the function~$g_n$ we have
\begin{equation}
\mdiv(g_n)=-m\infty^{(1)}-(m-1)\infty^{(2)}+\sum_{j=1}^{N-1}\maa_{N,j},
\label{35}
\end{equation}
where, as is clear, the zeros $\maa_{N,j}$ of the function $g_n$ can be chosen arbitrarily,
because the polynomials $q_{m,1},q_{m,2}$ are arbitrary. Next, the
function $g_n$ is meromorphic on~$\RS_2$ and the genus of the Riemann surface $\RS_2$ is zero,
and hence the function $g_n$ is completely defined by its divisor~\eqref{35} (of the zeros and poles).
As a~result, from~\eqref{35} we have the following explicit representation for the function $g_n$:
\begin{equation}
g_n(\zz)=
C_N\cdot \prod_{j=1}^{N-1} \bigl[\Phi(\zz)-\Phi(\maa_{N,j})\bigr]
\cdot \Phi(\zz)^{-m+1},
\quad C_N\neq0.
\label{35.2}
\end{equation}
Indeed, it is easily checked that the divisor of the zeros and poles of the
right-hand side of~\eqref{35.2} coincides with that of~\eqref{35}.
Below, in accordance with~\eqref{33}, we shall need to consider only the case
when all points~$\maa_{N,j}$ lie on the second sheet of the Riemann surface $\RS_2$,
$\maa_{N,j}=a^{(2)}_{N,j}\in\RS^{(2)}$. More precisely, the zeros $\maa_{N,j}$ should be
as follows: they should lie on the second list and be such that $\pi(\maa_{N,j})\in \myh{F}\setminus{E}$,
where $\myh{F}$ is the convex hull of~$F$.
In this case, it follows from~\eqref{35.2} that
\begin{equation}
g_n(z^{(2)})\Phi(z^{(2)})^{n+m+1}
=C_N\cdot \prod_{j=1}^{N-1}\bigl[\Phi(z^{(2)})-\Phi(a^{(2)}_{N,j})\bigr]
\cdot \Phi(z^{(2)})^{n+2}.
\label{35.3}
\end{equation}
We now consider  the product $g_n(\zz)\Phi(\zz)^{n+m+1}$. Using identities~\eqref{4} and~\eqref{4.3},
we write it as
\begin{align}
g_N(\zz)\Phi(\zz)^{n+m+1}
&=C_N\cdot \prod_{j=1}^{N-1}\bigl[\Phi(\zz)-\Phi(\maa_{N,j})\bigr]
\cdot \Phi(\zz)^{-m+1}\Phi(\zz)^{n+m+1}\notag\\
&=\myt{C}_N\cdot\prod_{j=1}^{N-1}\frac{z-a_{N,j}}{1-\Phi(\zz)\Phi(\maa_{N,j})}
\cdot\Phi(\zz)^{N+n+1},
\label{36}
\end{align}
where $\myt{C}_N\neq0$ and it is assumed that all
$\maa_{N,j}\neq\infty^{(1)},\infty^{(2)}$. In accordance with~\eqref{33},
we shall require representation~\eqref{36} only in the case when $\zz=z^{(2)}$
and all $\maa_{N,j}=a^{(2)}_{N,j}$. In this setting, we have by~\eqref{36}
\begin{equation}
g_N(z^{(2)})\Phi(z^{(2)})^{n+m+1}
=\myt{C}_N\prod_{j=1}^{N-1}\frac{z-a_{N,j}}{1-\Phi(z^{(2)})\Phi(a^{(2)}_{N,j})}
\cdot\Phi(z^{(2)})^{N+m+1}.
\label{37}
\end{equation}
Since  $\Phi(z^{(2)})=1/\pfi(z)$ for all $z\in D$, the last relation can be written as
\begin{equation}
g_N(z^{(2)})\Phi(z^{(2)})^{n+m+1}
=C_3(N)\prod_{j=1}^{N-1}\frac{z-a_{N,j}}{1-\pfi(z)\pfi(a_{N,j})}
\cdot\frac1{\pfi^{n+2}(z)}.
\label{38}
\end{equation}
Using \eqref{38}, the orthogonality relation~\eqref{33} can be put in the form
\begin{equation}
\int_F Q_{n,2}(x)\prod_{j=1}^{N-1}\frac{x-a_{N,j}}{1-\pfi(x)\pfi(a_{N,j})}
\cdot\frac{\pfi'(x)}{\pfi^{n+2}(x)}\,d\sigma(x)=0,
\label{39}
\end{equation}
where the number $N\leq{n}$ is arbitrary and all points $a_{N,j}$ lie in~$D$. From~\eqref{39},
it follows that $\mdeg{Q_{n,2}}=n$, all zeros of the polynomial $Q_{n,2}$ lie on~$\myh{F}$
(which is the convex hull of the compact set~$F$); besides, the gap with number $(p-1)$ between the intervals
$[c_j,d_j]$, $j=1,2,\dots,p$, may contain at most $p-1$ zeros of this polynomial.
The orthogonality relations~\eqref{39}, which are defined for an arbitrary
$N\leq{n}$ and arbitrary points $a_{N,j}\in\myh{F}\setminus{E}$, will underlie our further analysis.

\subsection{}\label{s3s3}
As usual, when applying\footnote{Note that, under the hypotheses of Theorem~\ref{th2},
the $\GRS$-method is much easier to deal with, because an $S$-compact set~$F$
is a~finite union of intervals of the real line and $\sigma$~is a~positive measure on~$F$; cf.\
\cite{Sta97b},~\cite{GoRa87},~\cite{Rak18}.}
the $\GRS$-method, we assume that
\begin{equation}
\frac1n\chi(Q_{n,2})\not\to\lambda=\lambda_F
\label{39.2}
\end{equation}
as $n\to\infty$.
We shall arrive at a~contradiction by using the orthogonality relations~\eqref{39} and
condition~\eqref{39.2}.

The weak compactness of the space of measures $M_1(\myh{F})$ shows that
\begin{equation}
\frac1n\chi(Q_{n,2})\to\mu\neq\lambda,
\quad n\in\Lambda,\quad n\to\infty
\label{40}
\end{equation}
for some infinite
subsequence $\Lambda\subset\NN$;
besides, $S(\mu)\subset{F}$,
$\mu\in M_1(F)$, $\mu(1)=1$ by the above properties of the polynomial $Q_{n,2}$.
We claim that relation~\eqref{40} and the orthogonality relation~\eqref{39}
contradict  each other.

Setting
$$
\myt{V}^\mu(z):=\int_F\log\frac1{|1-\pfi(z)\pfi(t)|}\,d\mu(t),
$$
we have
$$
P^\mu(z)=2V^\mu(z)-\myt{V}^\mu(z).
$$
Since $\mu\neq\lambda$, it follows that, for $z\in S(\mu)\subset F$,
\begin{equation}
P^\mu(z)+\psi(z)\not\equiv m_0
:=\min_{z\in F}\bigl(P^\mu(z)+\psi(z)\bigr)=P^\mu(x_0)+\psi(x_0),
\label{41}
\end{equation}
where $x_0\in F$. Hence there exists a~point $x_1\in S(\mu)$, $x_1\neq
x_0$, and a~number $\eps>0$ such that
\begin{equation}
P^\mu(x_1)+\psi(x_1)=m_1>m_0+\eps.
\label{41.2}
\end{equation}
Further, since the function $\psi(z)$ is harmonic and
the potential $P^\mu$ is lower semi-\allowbreak continuous, the same
inequality~\eqref{41.2} holds in some $\delta$-neighbourhood
$U_\delta(x_1):=(x_1-\delta,x_1+\delta)\not\ni x_0$, $\delta>0$, of the point~$x_1$.
We have $x_1\in S(\mu)$, and so $\mu(U_\delta(x_1))>0$. Hence, for all sufficiently large
$n\geq n_0$, $n\in\Lambda$, there exists a~polynomial
$p_n(z)=(z-\zeta_{n,1})(z-\zeta_{n,2})$ such that $\zeta_{n,1},\zeta_{n,2}\in
U_\delta(x_1)$ and $p_n$ divides the polynomial $Q_{n,2}$; that is,  $Q_{n,2}/p_n\in\PP_{n-2}$.
We set
\begin{equation}
\myt{Q}_n(z):=\frac{Q_{n,2}(z)}{p_n(z)}=\prod_{j=1}^{n-2}(z-x_{n,j}).
\label{41.3}
\end{equation}
We may assume in what follows that, for $n\in\Lambda$, all zeros of the polynomial
$Q_{n,2}$ lie in the set $\myh{F}\setminus{E}$.
Indeed, there is at most one gap between the intervals $[c_j,d_j]$ that may
contain the interval~$E$, in each gap lying at most one zero of the polynomial
$Q_{n,2}$. If some zero of the polynomial $Q_{n,2}$ lies on the interval~$E$, then in
definition~\eqref{41.3} of the polynomial $\myt{Q}_n$ one should replace
the corresponding factor ($(z-x_{n,j_0})$, say) by the factor $(z-\myt{x}_{n,j_0})$,
where the point $\myt{x}_{n,j_0}$ still lies in the (open) gap, but it is
not lying in~$E$ anymore.

Now in the orthogonality relation~\eqref{39} we put $N=n-1$  and
take the zeros $x_{n,j}$ of the polynomial $\myt{Q}_n$ as points $a_{N,j}$
(with the possible correction mentioned above),
relation \eqref{39} assuming the form
\begin{align}
0&=
\int_{F\setminus{U_\delta(x_1)}}
\frac{Q_{n,2}^2(x)}{p_n(x)}\prod_{j=1}^{n-2}\frac1{1-\pfi(x)\pfi(x_{n,j})}
\cdot\frac{\pfi'(x)}{\pfi^{n+2}(x)}\,d\sigma(x)\notag\\
&+
\int_{\myo{U}_\delta(x_1)}
\frac{Q_{n,2}^2(x)}{p_n(x)}\prod_{j=1}^{n-2}\frac1{1-\pfi(x)\pfi(x_{n,j})}
\cdot\frac{\pfi'(x)}{\pfi^{n+2}(x)}\,d\sigma(x).
\label{42}
\end{align}
We denote by $I_{n,1}$ and $I_{n,2}$, respectively, the first and second integrals in~\eqref{42}.
Since the integrand
in $I_{n,1}$ has constant sign for $x\in F\setminus{U_\delta(x_1)}$, we have
\begin{align}
|I_{n,1}|
&=
\int_{F\setminus{U_\delta(x_1)}}\biggl|
\frac{Q_{n,2}^2(x)}{p_n(x)}\prod_{j=1}^{n-2}\frac1{1-\pfi(x)\pfi(x_{n,j})}
\cdot\frac{\pfi'(x)}{\pfi^{n+2}(x)}
\biggr|\,d\sigma(x)\notag\\
&=\int_{F\setminus{U_\delta(x_1)}}
|Q_{n,2}(x)|\prod_{j=1}^{n-2}\biggl|\frac{x-x_{n,j})}
{1-\pfi(x)\pfi(x_{n,j})}\biggr|
\cdot\frac{\pfi'(x)}{\pfi^{n+2}(x)}\,d\sigma(x).
\label{43}
\end{align}
A similar analysis (see Lemma~7 of~\cite{GoRa87})  with the use of standard machinery of the  logarithmic potential theory
shows that
\begin{equation}
\lim_{\substack{n\to\infty\\n\in\Lambda}} |I_{n,1}|^{1/n}
=\exp\biggl\{
-\min_{x\in F\setminus{U_\delta(x_1)}}\bigl(P^\mu(x)+\psi(x)
\bigr)\biggr\}
=e^{-m_0}.
\label{44}
\end{equation}
We give a proof of~\eqref{44} for completeness  (cf.~Lemma~7 of~\cite{GoRa87}).

Indeed,
\begin{equation}
-\frac1n\sum_{j=1}^{n-2}\log|1-\pfi(x)\pfi(x_{n,j})|
\to\int_F\log\frac1{|1-\pfi(x)\pfi(t)|}\,d\mu(t)=\myt{V}^\mu(x)
\label{71}
\end{equation}
as $n\to\infty$ uniformly in~$x\in F$. Hence,
\begin{align}
\min_{x\in F}\biggl\{-\frac1n
\log\biggl(
|Q_{n,2}(x)|\prod_{j=1}^{n-2}\biggl|\frac{x-x_{n,j}}
{1-\pfi(x)\pfi(x_{n,j})}\biggr|\cdot\frac{\pfi'(x)}{\pfi^{n+2}(x)}
\biggr)\biggr\} \nonumber \\
\to\min_{x\in F} \bigl\{P^\mu(x)+\psi(x)\bigr\}
\label{72}
\end{align}
as $n\to\infty$. As a~result, we have
\begin{align}
\max_{x\in F}\biggl\{|Q_{n,2}(x)|\prod_{j=1}^{n-2}\biggl|\frac{x-x_{n,j}}
{1-\pfi(x)\pfi(x_{n,j})}\biggr|\cdot\frac{\pfi'(x)}{\pfi^{n+2}(x)}
\biggr\}^{1/n} \nonumber \\
\to\exp\bigl\{-\min_{x\in F}\bigl[P^\mu(x)+\psi(x)\bigr]\bigr\}
\label{73}
\end{align}
as $n\to\infty$, proving thereby the upper estimate
\begin{equation*}
\varlimsup_{\substack{n\to\infty\\n\in\Lambda}} |I_{n,1}|^{1/n}
\leq e^{-m_0}.
\end{equation*}

Let us now prove the corresponding lower estimate.
The potential $P^\mu$ is weakly continuous, and hence the function $P^\mu+\psi$
is approximately continuous with respect to the Lebesgue measure on the compact set~$F$.
Consequently, for any $\eps>0$, the set
$$
e=\{x\in F: (P^\mu+\psi)(x)<m_0+\eps\}
$$
has positive Lebesgue measure. From our assumptions we have
$$
-\frac1n\log
\biggl\{|Q_{n,2}(x)|\prod_{j=1}^{n-2}\biggl|\frac{x-x_{n,j}}
{1-\pfi(x)\pfi(x_{n,j})}\biggr|\cdot\frac{\pfi'(x)}{\pfi^{n+2}(x)}
\biggr\}
\to (P^\mu+\psi)(x)
$$
as $n\to\infty$ with respect to the measure on~$F$. So, the measure of the set
$$
e_n:=\biggl\{x\in e:
-\frac1n\log
\biggl(|Q_{n,2}(x)|\prod_{j=1}^{n-2}\biggl|\frac{x-x_{n,j}}
{1-\pfi(x)\pfi(x_{n,j})}\biggr|\cdot\frac{\pfi'(x)}{\pfi^{n+2}(x)}
\biggr)
<m_0+\eps\biggr\}
$$
tends to the measure of~$e$  as $n\to\infty$. Hence
\begin{equation}
\varliminf_{\substack{n\to\infty\\n\in\Lambda}} |I_{n,1}|^{1/n}
\geq e^{-(m_0+\eps)}\lim_{\substack{n\to\infty\\n\in\Lambda}}
\(\int_{e_n}\pfi'(x)\,d\sigma(x)\)^{1/n}=e^{-(m_0+\eps)},
\label{75}
\end{equation}
the last equality in~\eqref{75} holding because $\sigma'(x)>0$ a.e.
on~$F$. The lower estimate
$$
\varliminf_{\substack{n\to\infty\\n\in\Lambda}} |I_{n,1}|^{1/n}
\geq e^{-m_0}
$$
follows from \eqref{75}, because $\eps>0$  is arbitrary.
This proves~\eqref{44}.

On the other hand, for the second integral $I_{n,2}$ we have the estimate
\begin{equation}
|I_{n,2}|\leq
\int_{\myo{U}_\delta(x_1)}
|Q_{n,2}(x)|\prod_{j=1}^{n-2}\biggl|\frac{x-x_{n,j}}
{1-\pfi(x)\pfi(x_{n,j})}\biggr|
\cdot\frac{\pfi'(x)}{\pfi^{n+2}(x)}\,d\sigma(x).
\label{45}
\end{equation}
Now an analysis similar to that above shows that
\begin{align}
\varlimsup_{\substack{n\to\infty\\n\in\Lambda}} |I_{n,2}|^{1/n}
\leq\exp\biggl\{-\min_{x\in \myo{U}_\delta(x_1)}
\bigl(P^\mu(x)+\psi(x)\bigr)
\biggr\}
\leq e^{-m_1}<e^{-(m_0+\eps)}.
\label{46}
\end{align}
But relations \eqref{44} and~\eqref{46} contradict the equality $I_{n,1}=-I_{n,2}$,
which is consequent on the orthogonality relations~\eqref{39}.

This proves Theorem~\ref{th2}.

\begin{remark}\label{rem1}
In a certain sense, the  above transformations
mean the change of the variable~$z$ by the variable $\zeta=\pfi(z)$.
For the case of the Riemann surface $w^2=z^2-1$ under consideration,
the key orthogonality relation~\eqref{39} can be derived
directly from relations~\eqref{24}, properties of functions of the second
kind, and identity~\eqref{4.2}, which is much faster.
However, in the present paper we chose a~different method of exposition,
because in a~more general setting,
when, for example,
$$
f_1(z):=\int_{-1}^1\frac{r(x)}{(z-x)}\frac{dx}{\sqrt{1-x^2}},
$$
where $r\in\CC(z)$ is an arbitrary complex rational function without poles and
zeros on~$E$, such a~simplification does not apply anymore, but the conclusions of Theorem~\ref{th2}
remain valid\footnote{As was already mentioned above, the author
intends to investigate this general case in a~separate paper; see \cite{Sue18}.}.
It is also worth pointing out the role of the Riemann surface of the function $w^2=z^2-1$
in our analysis, because we also intend to extend both the results from the present paper and
those announced in~\cite{Sue18}
to the hyperelliptic setting, when, instead of the Riemann surface of the function $w^2=z^2-1$
of genus $g=0$, use is made of the Riemann surface of the function $w^2=(z-e_1)\dots(z-e_4)$
of genus  $g=1$.
A~generalization of the results obtained here to the elliptic case (of
course, if such an extension will come to being) will be of the utmost importance
in assessing the
potency of the method proposed here when investigating the general case of a~pair of functions
$f_1,f_2\in\sA^\circ(\Sigma)$. Of course, in this general case the problem of the
{\it formula for strong asymptotics} for Pad\'e polynomials valid
for an {\it arbitrary} function~$f$ from the class $\sA^\circ(\Sigma)$
will have a~great value;
see \cite{Nut86}, \cite{MaRaSu12},~\cite{ApYa15} in this respect.
\end{remark}

\begin{remark}\label{rem2}
It is well known (see~\cite{Apt99}, and also \cite{ApLy10} and~\cite{LoVa16}) that,
for a~pair of functions $f_1,f_2$ forming a~Nikishin system, the
support of the equilibrium measure~$\lambda$ in the diagonal case (which is considered here)
coincides with the entire compact set~$F$; this means that only the case of identical equality in relations~\eqref{14}
is possible. So far, this fact has not yet been proved within the framework of the approach proposed here.
\end{remark}

\begin{remark}\label{rem3}
It is worth pointing out that the approach proposed in the present paper  stems, to some extent,
from the analysis of numerical experiments of \cite{KoIkSu15}--\cite{KoIkSu16}.
\end{remark}

\clearpage
\newpage

\clearpage
\newpage

\clearpage


\end{document}